\documentclass[%
 reprint,
 amsmath,amssymb,
 aps,
pre,
]{revtex4-2}

\usepackage{dcolumn}
\usepackage{bm}
\usepackage{braket}
\usepackage{mathtools}
\usepackage{tensor}
\usepackage{graphicx} 
\allowdisplaybreaks

\usepackage{mleftright} 
\usepackage[super]{nth} 

\newcommand{\nx}{N_X}
\newcommand{\ny}{N_Y}
\newcommand{\matC}{\mathbf{C}}
\newcommand{\Csq}{\mathbf{C}^\top\mathbf{C}}
\newcommand{\Csqtranspose}{\mathbf{C}\mathbf{C}^\top}
\newcommand{\matH}{\mathbf{H}}

\usepackage[T2A,T1]{fontenc}
\usepackage[utf8]{inputenc}
\usepackage[russian,english]{babel}

\newcommand{\nn}{\nonumber}

\begin{document}

\preprint{APS/123-QED}

\title{Distribution of singular values in large sample cross-covariance matrices}

\author{Arabind Swain}
 \affiliation{Department of Physics, Emory University, Atlanta, GA 30322, USA}
\author{Sean Alexander Ridout}%
 \affiliation{Department of Physics, Emory University, Atlanta, GA 30322, USA}
 \affiliation{Initiative in Theory and Modeling of Living Systems, Atlanta, GA 30322, USA}
\author{Ilya Nemenman}
\affiliation{Department of Physics, Emory University, Atlanta, GA 30322, USA}
\affiliation{Department of Biology, Emory University, Atlanta, GA 30322, USA}
\affiliation{Initiative in Theory and Modeling of Living Systems, Atlanta, GA 30322, USA}

\date{\today}
\begin{abstract}
  For two high-dimensional datasets $X$ and $Y$, with dimensionalities $N_X$ and $N_Y$ of order of the number of samples $T$, estimates of their cross-covariance will have large fluctuations. These sampling fluctuations can be studied by analyzing the case of uncorrelated $X$ and $Y$, samples of which comprise large matrices ${\mathbf X}$ and ${\mathbf Y}$ with Gaussian i.i.d.\ entries and dimensions $T\times N_X$ and $T\times N_Y$, respectively. For this problem, we derive the probability distribution of the singular values of $\mathbf{X}^\top \mathbf{Y}$ in different parameter regimes. This extends the Marchenko–Pastur result for the distribution of eigenvalues of empirical sample covariance matrices to singular values of empirical cross-covariances. We analyze these results in a variety of limits, arguing that in many cases signals may be detected even if one or both dataset are of dimensionality greater than the number of samples, where methods based on whitening of the cross-covariance cannot be used. Our results will help to establish statistical significance of cross-correlations in many data-science applications.
\end{abstract}
\maketitle

\section{Introduction}
Many data-science applications require detecting correlations between two variables $X$ and $Y$ of dimensions $N_X$ and $N_Y$, respectively, with $N_X, N_Y \gg 1$. When these variables are sampled $T$ times, with $T \sim N_X, N_Y$, resulting in data matrices $\mathbf{X}\in R^{T\times N_X}$ and $\mathbf{Y}\in R^{T\times N_Y}$, respectively, sampling fluctuations can produce spurious correlations, even when $X$ and $Y$ are truly uncorrelated. Characterizing these sampling-induced correlations is essential before isolating genuine signals in real datasets.

Marchenko and Pastur famously analyzed similar correlations in sample self-covariance matrices  \cite{marchenko1967распределение} using now-classic methods of  Random Matrix Theory (RMT) \cite{Potters_Bouchaud_2020}. They derived the spectra of so-called Wishart matrices $\frac{1}{T}{\mathbf X}^\top{\mathbf X}$, where all entries of ${\mathbf X}$ are i.i.d.\ normal random variables. For $T > N_X, N_Y$, later work generalized these results to cross-correlations of large-dimensional whitened variables \cite{Bouchaud2007,Potters2023,firoozye2023canonical}, where {\em whitening} denotes linearly transforming data to zero mean and unit covariance matrix, so that there are no correlations remaining  within the transformed  $X$ and $Y$ individually; this parallels Canonical Correlation Analysis (CCA)~\cite{Hotteling}. However, when $T<N_X, N_Y$, whitening is non-trivial since the $X$--$X$ and $Y$--$Y$ self-covariance matrices cannot be inverted. In this case, Partial Least Squares (PLS)~\cite{wold1966estimation}, which deals with cross-correlations between unwhitened data, becomes essential. This regime is common in  many cases, where the number of samples is limited (see, e.g., \cite{YOU2024108434,LêCaoRossouwRobertGraniéBesse2008}). 

While the whitened case is well understood, to our knowledge, no similar explicit understanding exists for the unwhitened cross-covariance between $X$ and $Y$ for arbitrary values of $T, N_X, N_Y\gg1$. That is, even though many relevant  RMT results are known, no explicit expressions for the singular value spectra of cross-covariance have been written down, and the limits of these expressions for different regimes relevant for data analysis have not been explored. More specifically, in RMT, the spectrum of a random matrix $\mathbf{A}$ is usually obtained from its Stieltjes transform $\mathfrak{g}_{\mathbf{A}}(z)$ (see below for details), and several publications  obtained expressions for algebraic equations that can be solved to find the  Stieltjes transform of a product of two Wishart matrices~\cite{Philiplinear, Rocks2022,PhysRevE2010ref}, which, as we explain below, is a useful model for understanding spectra of cross-covariance of two datasets. In fact, similar results exist for products of arbitrarily many Wishart matrices, i.e.,  $(\mathbf{X}_1 \dots \mathbf{X}_M )(\mathbf{X}_1 \dots \mathbf{X}_M )^\top$, for both complex and real elements of ${\mathbf X}_m$~\cite{Forrester_2014,dupic2014spectraldensityproductswishart}. Some results have even been obtained for random matrices of the form $\sigma(\mathbf {WX})\sigma(\mathbf {WX})^\top$, where $\sigma(\cdot)$ is a nonlinear function, which arise in the context  of large neural networks~\cite{Benigni_2001,NIPS2017_Pennington}. However, none of these previous publications explicitly study consequences of their RMT calculations in the context of cross-covariance-based data analysis.

In this paper, we apply existing RMT methods to explicitly calculate and analyze  singular value spectra of unwhitened sample cross-covariance matrices, for uncorrelated Gaussian i.i.d.\ data and arbitrary relations among $T$, $N_X$, and $N_Y$. In particular, our results suggest that correlations between the variables  may be detectable even if the dimensionality of one or both variables is larger than  the number of samples, where CCA-like methods, which require inversion of marginal covariances ${\mathbf X}^\top{\mathbf X}$ and ${\mathbf Y}^\top{\mathbf Y}$  cannot be used. We hope that our results can be used to improve understanding of statistical significance of cross-correlations in  data science applications.

\section{Model and methods}
We consider $T$ samples of random variables $X$ and $Y$ combined into matrices $\mathbf X$ and $\mathbf Y$, with dimensions $T \times \nx$ and $T \times \ny$, respectively.  The entries of $\mathbf{X}$ and $\mathbf{Y}$ are i.i.d.\ Gaussian random variables with zero mean and variances $\sigma_X^2$ and $\sigma_Y^2$ respectively,
\begin{align}
   & X_{t\mu} \sim \mathcal N(0,\sigma_X^2)\,,\quad Y_{t \nu} \sim \mathcal N(0,\sigma_Y^2)\,,\label{eq:vars}\\
    &t=1,\ldots,T,\; \mu=1,\ldots,\nx,\; \nu=1,\dots,\ny.
\end{align}
In this model there are no true correlations between $X$ and $Y$, so the sample estimates of the correlation, computed from $\mathbf{X}$ and $\mathbf{Y}$ will be small when $N_X/T, N_Y/T$ are small.

We define normalized matrices as
\begin{equation}
    \widetilde{\mathbf X} = \frac{\mathbf{X}}{\sigma_X},\quad\quad
    \widetilde{\mathbf Y} = \frac{\mathbf{Y}}{\sigma_Y}\label{eqn:scaled_XY}\,.
\end{equation}
For $T\gg1$, each column in these matrices has variance of nearly one. Note that, in typical applications, the variance $\sigma^2_{X}$ of $X$ and the variance $\sigma^2_{Y}$ of $Y$ would be estimated from samples as well, and the estimates might be different from their true value. Here we disregard this distinction, as in \cite{Philiplinear}, arguing that sampling fluctuations in estimating scalar parameters are negligible compared to sampling effects on the thermodynamically many singular values.

The normalized empirical cross-covariance matrix (NECCM) $\mathbf C $ is then
\begin{equation}\label{eq:NECCM1}
    \matC=\frac{1}{T}\widetilde{\mathbf Y}^\top\widetilde{\mathbf X},
\end{equation}
which has dimensions $\ny \times \nx$.  
If $\nx \neq \ny$, this matrix is not square, but it obviously has the same nonzero singular values as its transpose. Without loss of generality, in all calculations, we take $\nx\leq \ny$.

We want to calculate the distribution of these singular values. To utilize RMT methods, most of which only work for square symmetric matrices, we focus instead on eigenvalues of 
\begin{equation}\label{CC1}
    \Csq=\frac{1}{T^2}\widetilde{\mathbf X}^\top\widetilde{\mathbf Y}\widetilde{\mathbf Y}^\top\widetilde{\mathbf X}.
\end{equation}
Nonzero eigenvalues of $\Csq$, which we denote as $\lambda$, are the same as nonzero eigenvalues of $\Csqtranspose$, and their distribution is related to the distribution of nonzero singular values of $\matC$, denoted as $\gamma$, via
\begin{equation}
    \rho_{C}(\gamma)=2\sqrt{\lambda}\rho_{C^\top C}(\lambda),\quad \gamma=\sqrt{\lambda}. \label{eq:jac}
\end{equation}

The matrices $\mathbf{C}^\top\mathbf{C}$ and $\mathbf{C}\mathbf{C^\top}$ have the same nonzero eigenvalues, with density denoted by $\tilde{\rho}(\lambda)$. The distribution of eigenvalues of $\Csq$ will further contain a $\delta$-function at zero consisting of $\nx-T$ zero eigenvalues if $T\leq \nx$ (recall that we assume $N_X \leq N_Y)$. Thus,
\begin{equation}
    \rho_{\Csq}(\lambda)=\frac{\min (N_X,T)}{N_X}\tilde{\rho}(\lambda)+\left(1-\frac{\min (N_X,T)}{N_X}\right)\delta(\lambda), \label{eq:rho_1}
\end{equation}
The distribution of eigenvalues of $\Csqtranspose$ will contain $N_Y-N_X$ additional zero eigenvalues. Thus,
\begin{equation}
    \rho_{\Csqtranspose}(\lambda)=\frac{\min (N_X,T)}{N_Y}\tilde{\rho}(\lambda)+\left(1-\frac{\min (N_X,T)}{N_Y}\right)\delta(\lambda), \label{eq:rho_2}
\end{equation}

To explore the problem in different regimes, we define:
\begin{align}
    q_X \equiv \nx/T, \; q_Y\equiv \ny/T,\;p_X\equiv1/q_X, \; p_Y\equiv1/q_Y.
\end{align}
Our RMT results for the spectrum will hold in the limit $N_X, N_Y, T \to \infty$, with $p_X$ and $p_Y$ held fixed.

\emph{Eigenvalue density.} We compute the eigenvalue density of the square of the NECCM, Eq.~\eqref{CC1}, by computing its Stieltjes transform, as is the standard approach \cite{Potters_Bouchaud_2020}. The Stieltjes transform of an $N \times N$ matrix $\mathbf{A}$, with eigenvalues $\lambda_1, \dots, \lambda_N$, is defined as
\begin{equation}\label{Stieltjes}
    g_{A,N}(z)=N^{-1}\text{Tr}(z\mathbf I-\mathbf A)^{-1}=N^{-1}\sum_{i=1}^N\frac{1}{z-\mathbf{\lambda}_i},
\end{equation}
 where $z$ is a complex number, which is restricted to either positive or negative imaginary part so as to be defined away from all the (real) eigenvalues of $\mathbf{A}$. We denote the large-$N$ limit of $g_{A,N}$ by $\mathfrak g_\mathbf{A}$~\cite{Potters_Bouchaud_2020}, $\mathfrak g_{\mathbf{A}}(z) = \lim_{N\to \infty} \mathbb{E}[\mathfrak g_{\mathbf{A},N} (z)]$. The  eigenvalue density is obtained from the Sokhotski–Plemelj formula
\begin{align}\label{eq:sokhotski-plemelj}
    \rho_A(\lambda)=\lim_{\eta\rightarrow 0^+}\frac{1}{\pi}\Im \mathfrak g_{\mathbf A}(\lambda - i\eta)\,,
\end{align}
where $\Im$ denotes the imaginary part. We use a series of relatively common random matrix operations to obtain the Stieltjes transform of the square of NECCM, in the limit where $\nx, \ny, T \to\infty$ with $p_X$ and $p_Y$ held fixed. These steps are outlined in the {\em Appendix}.

In general, we find that $\tilde{\rho}(\lambda)$ is nonzero over some finite interval $(\lambda_-, \lambda_+)$. The corresponding values for nonzero singular values of the NECCM are denoted by $\gamma_\pm$. As the imaginary part of the Stieltjes transform gives us the eigenvalue density of the square of the NECCM, $\lambda_{\pm}$ can be found by solving a discriminant equation associated with the algebraic equation for the Stieltjes transform (see {\em Appendix}). Analytical expression for these boundaries for the cross-covariance spectrum of pure uncorrelated  noise are one of the central results of this paper. 

\emph{Numerical simulations.} We confirm our results by simulating the model, Eq.~(\ref{eq:vars}), numerically. Although the eigenvalue density is expected to be self-averaging, and thus our calculations for $\rho(\gamma)$ will be exact for SVD of an \emph{individual} matrix for sufficiently large $T$, making $T$ very large substantially increases the computational costs. Thus, we simulate matrices with $T=1000$, and more precisely test our predictions by averaging over $500$ independent realizations.

\section{Equation for Stieltjes transform and singular value density bounds}
We calculate the density of eigenvalues of the square of NECCM in 3 cases, covering all  possible relationships between $T, \nx, \ny$: (1) $T>\nx,\ \ny$, (2) $\ny\geq T\geq \nx$, and (3) $T< \nx,\, \ny$. For analyzing these different cases, we note that the square of the NECCM can be written as an $\nx \times \nx$ matrix $\Csq=\frac{1}{T^2}\widetilde{\mathbf X}^\top\widetilde{\mathbf Y}\widetilde{\mathbf Y}^\top\widetilde{\mathbf X}$ or an $\ny \times \ny$ matrix $\mathbf C\mathbf C^\top=\frac{1}{T^2}\widetilde{\mathbf Y}^\top\widetilde{\mathbf X}\widetilde{\mathbf X}^\top\widetilde{\mathbf Y}$. Both of these matrices will have the same nonzero eigenvalues, as indicated in Eqs.~(\ref{eq:rho_1}, \ref{eq:rho_2}). Similarly, the $T \times T$ matrix  $\matH =\frac{1}{T^2}\widetilde{\mathbf X}\widetilde{\mathbf X}^\top\widetilde{\mathbf Y} \widetilde{\mathbf Y}^\top$ will have the same nonzero eigenvalues, i.e.,
\begin{equation}
    \rho_{\mathbf{H}}(\lambda)=\frac{\min (N_X,T)}{T}\tilde{\rho}(\lambda)+\left(1-\frac{\min (N_X,T)}{T}\right)\delta(\lambda). \label{eq:rho_H}
\end{equation}
Through Eqs.~(\ref{eq:rho_1}, \ref{eq:rho_2}, \ref{eq:rho_H}), all Stieltjes transforms can be related to the Stieltjes transform $h_T(z)\equiv g_{ \matH, T}(z)$ of $\mathbf{H}$, giving
 \begin{align}
     g_{\Csq, \nx}(z)&= p_X h_T(z)+\left(1-p_X\right)\delta(z),\label{eq:finiteN_Stieltjes} \\ 
    g_{\Csqtranspose, \ny}(z)&= p_Y h_T(z)+\left(1-p_Y\right)\delta(z).
\end{align}

An RMT calculation (Appendix \ref{append:sec1}) then shows that $\mathfrak{h}(z) \equiv \lim_{T\to\infty} h_T(z)$ satisfies a cubic equation 
\begin{equation}\label{eq:cubicstilt}
    a \mathfrak h^3+  b \mathfrak h^2+ c \mathfrak h+  d=0,
 \end{equation}
where
\begin{align}
    &a=z^2p_Xp_Y,\\
    &b=z\left(p_Y(1-p_X)+p_X(1-p_Y)\right),\\
    &c=\left((1-p_X)(1-p_Y)-z p_Xp_Y\right),\\
    &d=p_Xp_Y.
\end{align}
Thus, solving Eq.~(\ref{eq:cubicstilt}), and then using Eq.~(\ref{eq:finiteN_Stieltjes}), gives the eigenvalue density of $\Csq$, which can be used to compute the density of the nonzero singular values of the cross-covariance using Eq.~(\ref{eq:jac}).

\begin{figure*}
     \centering
     \includegraphics[width=\textwidth]{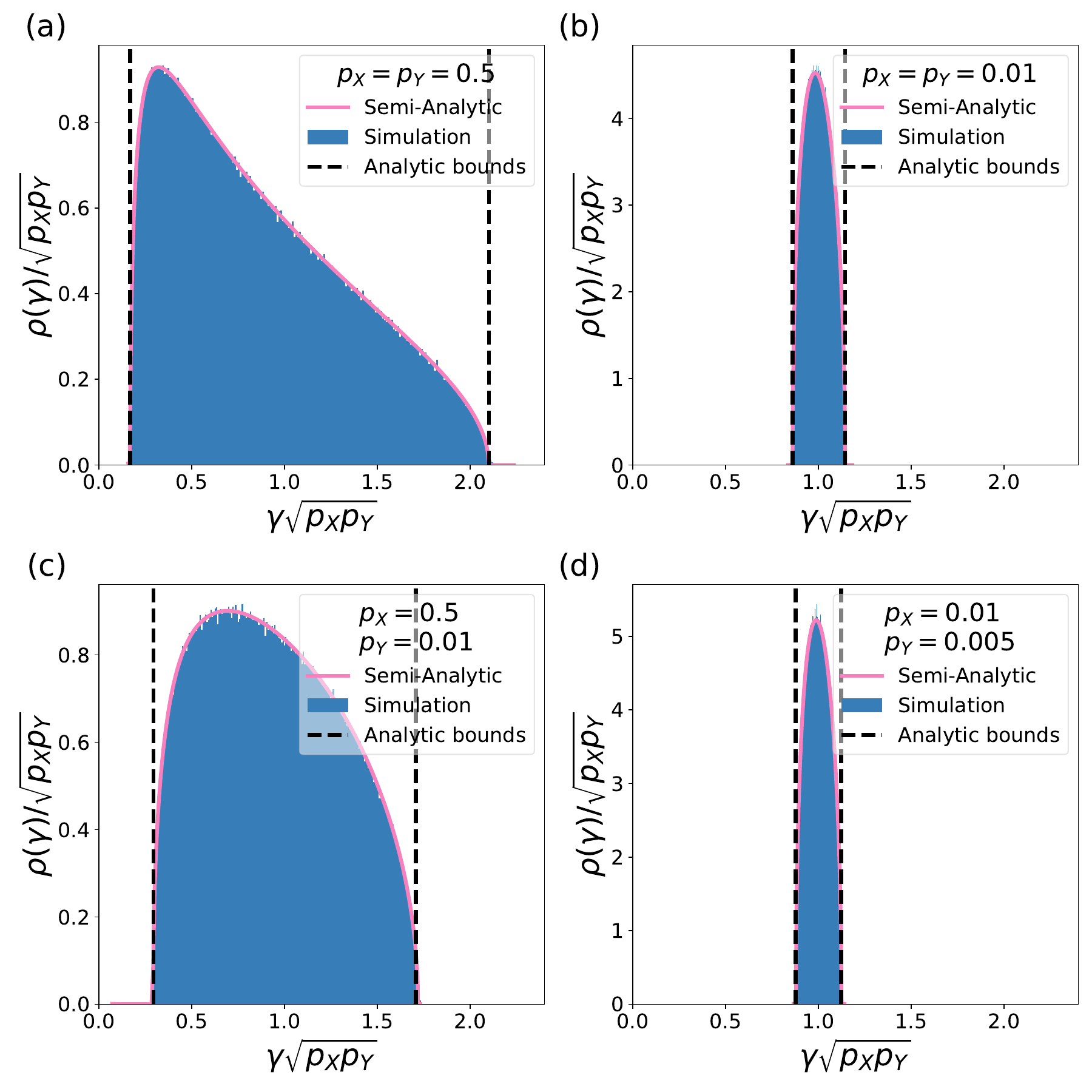}
     \caption{Distribution of nonzero singular values for $T < \nx ,\ \ny$. \textbf{(a)} $p_X=p_Y=0.5$ with analytic bounds given by Eq.~(\ref{eq:edge}), \textbf{(b)} $p_X=p_Y=0.01$ with analytic bounds given by Eq.~(\ref{eq:edge1}), \textbf{(c)} $p_X=0.5$, $p_Y=0.01$ with analytic bounds given by Eq.~(\ref{eq:edge2}), and \textbf{(d)} $p_X=0.01$, $p_Y=0.05$ with analytic bounds given by Eq.~(\ref{eq:edge3}). The blue bars are the histograms of the simulated data. The magenta curves are computed from the numerical solution of the exact cubic equation for the Stieltjes transform. The black dashed lines show bounds of the nonzero part of the density in simplifying limits, evaluated analytically. Here, $T=1000$, and the the numerical simulation for spectrum consists of $500$ independent model realizations. We scale the singular values by $\sqrt{p_Xp_Y}$.  This places the midpoint of the distribution within a factor of a few from 1, and the range of the distributions between 1 and 10, for all parameters explored here. }
     \label{fig:undersampled}
\end{figure*}

\subsection{Spectrum of empirical cross covariance
matrix when $T < \nx ,\ \ny$}

The cubic polynomial given by Eq.~(\ref{eq:cubicstilt}) can be solved, numerically or analytically, for the imaginary part of $\mathfrak{h}$ at any parameter values. Taking its imaginary part then gives us the density of nonzero eigenvalues.

Here, we solve the equation numerically (which we refer to as the ``semi-analytic'' solution, since it solves numerically the analytical expression, Eq.~(\ref{eq:cubicstilt})), and study the spectrum for a variety of parameter regimes. The spectrum has compact support, showing a single band of eigenvalues with upper and lower bounds. The  bounds can be calculated by finding the condition under which the the discriminant of the cubic equation, Eq.~(\ref{eq:cubicstilt}), becomes zero. To get easily interpretable formulas for the bounds $\lambda_{\pm}$  (and hence $\gamma_{\pm}$), we take various simplifying limits where the discriminant equation for the cubic polynomial is exactly solvable. 

Firstly, consider the case where $p_X=p_Y <1$ (same-size data matrices, with $T < N_X, N_Y$). In this case, the bounds of the spectrum simplify to
\begin{equation}
\gamma_{\pm}=\sqrt{\frac{8p_X^2 + 20 p_X^3 - p_X^4 \pm p_X^{5/2} (8 + p_X)^{3/2}}{8p_X^4}}.\label{eq:edge}
\end{equation}
(The generalization to all values of $p_X$ is given alongside the derivation in the Appendix, cf. Eq.~(\ref{eq:xyequal}).) 
 
Assuming $p_X=p_Y \to 0$ (so that we are in the severely undersampled regime, where the number of samples is {\em much} smaller than the number of dimensions in $X$ and $Y$), the edge values become \begin{equation}
\gamma_{\pm}\approx\frac{1}{p_X}(1\pm\sqrt{2p_X}).\label{eq:edge1}
\end{equation}

Secondly, we consider the case where one dataset is much higher-dimensional than the other, $p_Y \ll p_X \leq 1$. In this limit,
\begin{equation}
\gamma_{\pm}\approx\sqrt{\frac{1+p_X\pm2\sqrt{p_X}}{p_X p_Y}}.\label{eq:edge2}
\end{equation}

Finally, we can obtain simple results when both $p_X, p_Y \ll 1$, with $p_X/p_Y = O(1)$ (both $X$ and $Y$ are extremely undersampled, but unequal in dimensionality). In this limit, the bounds are \begin{equation}\gamma_{\pm}\approx\frac{1\pm\sqrt{p_Y+p_X}}{\sqrt{p_Yp_X}}.\label{eq:edge3}
\end{equation}

We see that, in all of these limits, magnitude of the singular values are roughly of the order of $\sqrt{\frac{1}{p_X p_Y}}=\sqrt{q_Xq_Y}$. This sets the typical scale of sampling noise  singular values at a given sample size $T$. The noise eigenvalues of $\tilde{\mathbf{X}}^\top\tilde{\mathbf{X}}/T$ and $\tilde{\mathbf{Y}}^\top\tilde{\mathbf{Y}}/T$ individually scale like $q_X$ and $q_Y$~\cite{marchenko1967распределение}. Thus, this scaling is plausible if each eigendirection is poorly-sampled enough that they can be found to correlate with each other by chance: evidently, since $N/T=O(1)$, this is the case.

Figure~\ref{fig:undersampled} compares our analytical results to numerical simulations for the density of singular values $\gamma$ of $\matC$. We scale the singular values by the  scale factor $\sqrt{\frac{1}{p_X p_Y}}$. We see that the semi-analytic solution for the density is in excellent agreement with our numerical results. Further, we see that the analytical solutions for the bounds, in appropriate limits, also agree well with simulations.

The simulations and the semi-analytic solutions also agree for other parameter values where simple analytic formulas for the bounds could not be evaluated exactly (see Appendix~\ref{append:sec1}).  

\subsection{Spectrum of empirical cross covariance matrix when $\ny \geq T\geq \nx$}
Solving for the roots of the cubic polynomial in  Eq.~(\ref{eq:cubicstilt}) and taking its imaginary part again gives us the density of nonzero eigenvalues. 

In this case, we can evaluate the bounds of the spectrum exactly in the limit $p_Y \ll 1 < p_X$. In this case, the bounds are \begin{equation}\gamma_{\pm}=\sqrt{\frac{1+p_X\pm2\sqrt{p_X}}{p_X p_Y}}.
\end{equation}
This limit is the same as in the case when $T \leq \nx ,\ \ny$, although the number of zero eigenvalues differs between the two cases.

\begin{figure}
    \centering
    \includegraphics[width=0.49\textwidth]{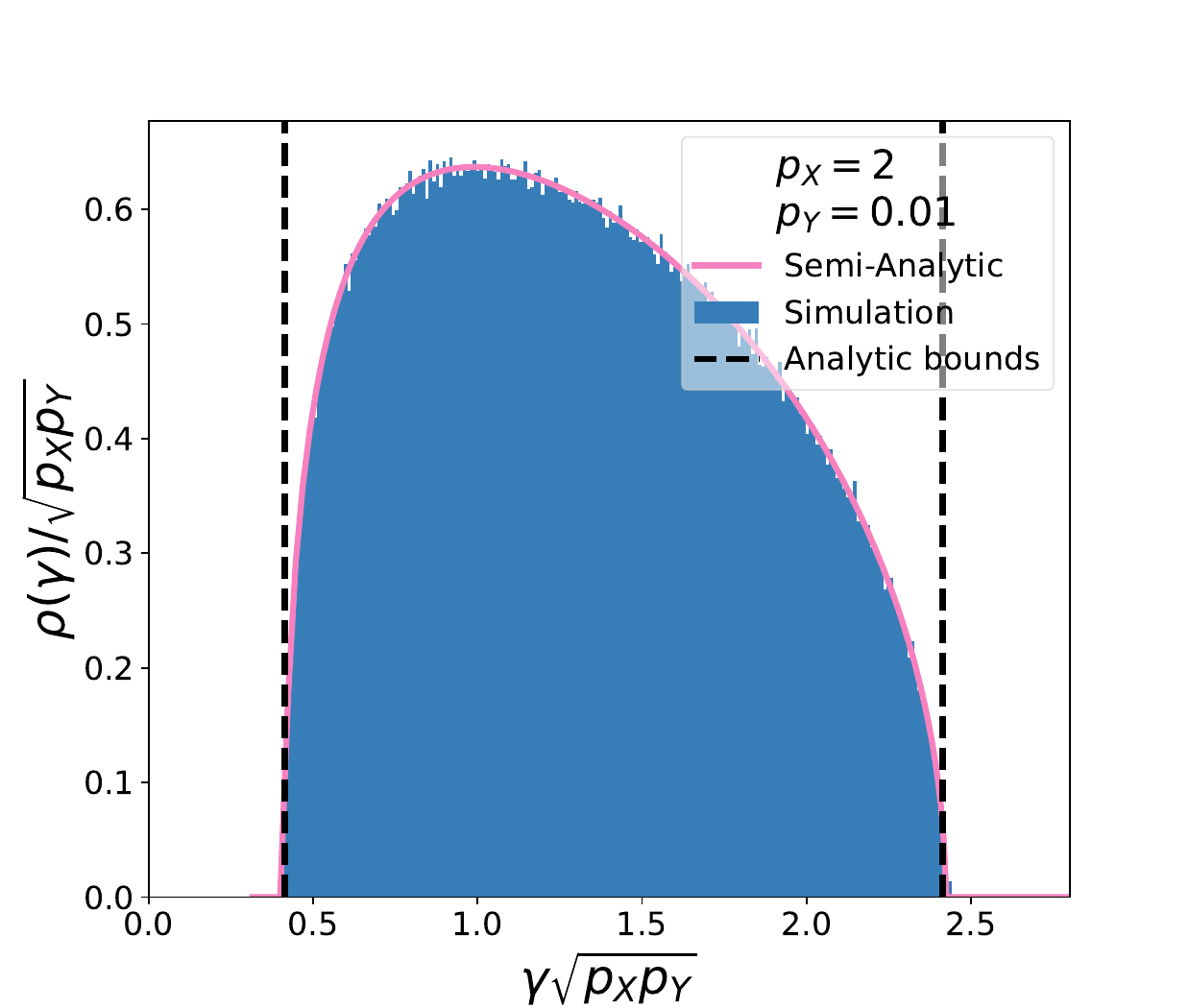}
    \caption{Distribution of nonzero singular values for $\ny \geq T\geq \nx$, specifically,  $p_X=2$, $p_Y=0.01$ with analytic bounds given by Eq.~(\ref{eq:edge2}). Plotting conventions are the same as in Fig.~\ref{fig:undersampled}. Here, again, $T=1000$, and the numerical simulation for spectrum consists of $500$ independent model realizations.}
    \label{fig:mixedlimit}
\end{figure}

Figure \ref{fig:mixedlimit} shows that the semi-analytic solution for the density, and the analytic solution for the bounds, match our numerical simulations in this case as well.

\subsection{Spectrum of empirical cross covariance matrix for $T>\nx,\ny$}\label{sec3c}

Solving for the roots of the cubic polynomial, Eq.~(\ref{eq:cubicstilt}), and taking its imaginary part again gives us the density of nonzero eigenvalues. We then obtain simplified formulas for $\gamma_{\pm}$ in  limiting cases.

Recall that the density of eigenvalues is nonzero when $\mathfrak{h}$ has an imaginary part. The boundaries of this region are identified by solutions $z$ of the discriminant of the cubic equation for $\mathfrak{h}$. For the case where $p_X=p_Y$, the discriminant is a $\nth{5}$-order polynomial with three zero solutions and two nonzero solutions $z_{\pm}$, where $z_{\pm}=\frac{8p_X^2 + 20 p_X^3 - p_X^4 \pm p_X^{5/2} (8 + p_X)^{3/2}}{8p_X^4}$. Now because $z_{-} <0$ and the squares of singular values are always positive, the upper bound of the nonzero density is $z_{+}$ but the lower bound is 0. Thus the bounds for the nonzero eigenvalue density are \begin{equation}\gamma_{+}=\sqrt{\frac{8p_X^2 + 20 p_X^3 - p_X^4 + p_X^{5/2} (8 + p_X)^{3/2}}{8p_X^4}}, \quad \gamma_{-}=0.\label{eq:edge4}
\end{equation}
In the limit $p_X\gg 1$ (extremely good sampling), this simplifies to $\gamma_{+}\approx \sqrt{\frac{3}{2p_X}}=\sqrt{\frac{3q_X}{2}}$. Thus, in this limit the scaling of the bounds agrees with those for the cross-correlations of whitened variables evaluated in Ref.~\cite{Bouchaud2007}, where $\gamma_{+}=2\sqrt{q_X}$, and $\gamma_{-}=0$. Note, however, that the exact value of the upper edge is different for the whitened cross-correlation matrices, because the self-covariances used for whitening also fluctuate. 

\begin{figure}
    \centering
    \includegraphics[width=0.49\textwidth]{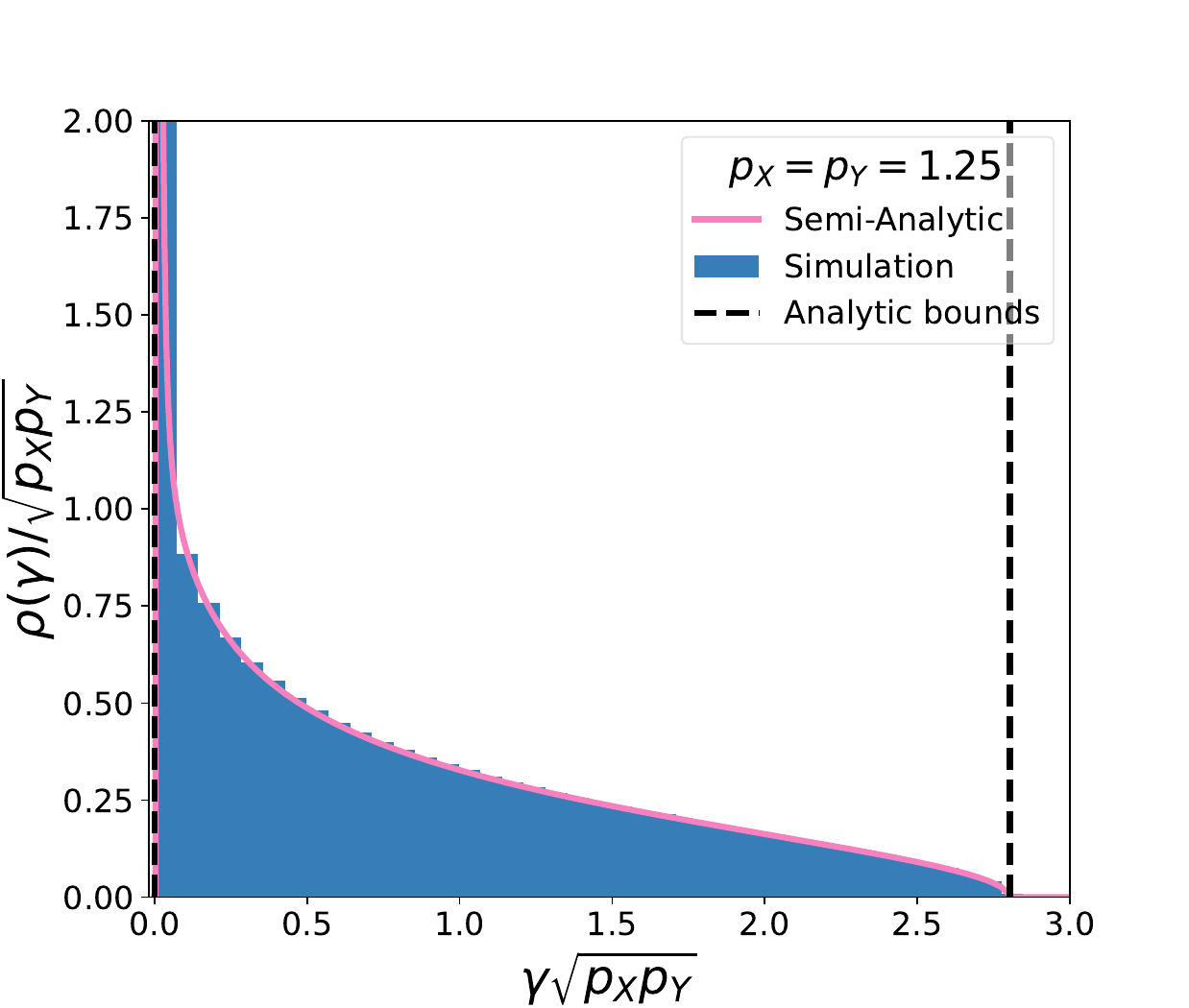}
    \caption{Distribution of nonzero singular values for $T>\nx,\ny$, specifically $p_X=p_Y=1.25$ with analytic bounds given by Eq.~(\ref{eq:edge4}). Plotting conventions are the same as in Fig.~\ref{fig:undersampled}. Here, again, $T=1000$, and numerical simulation for spectrum consists of $500$ independent model realizations.}
    \label{fig:classicallimit}
\end{figure}

Figure~\ref{fig:classicallimit} shows that these limiting formulas for the bounds, and the semi-analytic solution for the spectrum match numerical simulations.

\section{Discussion}
We have used random matrix theory to calculate the density of singular values of normalized cross-correlation matrices. Further, in simplifying limits, we were able to obtain simple, exact formulas for the bounds of the spectrum.

In all cases, the scale of the nonzero singular values is given roughly by $1/\sqrt{p_X p_Y} =\sqrt{\nx \ny} / T$. Thus, the noise, unsurprisingly, decreases as more samples are collected, relative to the dimensions of the two observed variables. More surprisingly, however, this calculation in fact suggests that the cross-covariance can sometimes be used to detect a signal which is not detectable from either the covariance of $X$ or that of $Y$ alone, as recently observed numerically~\cite{abdelaleem2024simultaneous}.

To see this, consider a naïve protocol for establish a correlation between high-dimensional $X$ and $Y$: we first search for a low-dimensional signal in $\mathbf{X}$ (e.g., using principle component analysis), then search for a low-dimensional signal in $\mathbf{Y}$, and finally correlate the low-dimensional signals. The bounds of the empirical covariance spectra of $\mathbf{X}$ and $\mathbf{Y}$ are of order $1/p_X$ and $1/p_Y$, respectively. Thus, a shared signal that has $O(1)$ magnitude in both $X$ and $Y$ will correspond to an outlier eigenvalue outside of the spectrum, and hence can be detected if  $T > \nx, \ny$. In particular, if $\ny > T > \nx$ (one variable is well sampled, and one variable is poorly sampled), the signal in $\mathbf{X}$ cannot be detected. Since the noise spectrum of $\matC$ depends on the geometric mean $\sqrt{p_X p_Y}$, however, the same signal may be detectable in $\matC$, if $X$ is sampled well enough to ``make up for'' the poor sampling of $Y$. Making this rough analysis precise requires a full calculation of the spectrum of a model with both a signal and noise, which we will present in a future work.

These results also suggest that a sufficiently strong signal can be detected even if $T < \nx, \ny$.

In the limit $T\gg \nx, \ny$, where the covariances of $X$ and $Y$ are both well sampled, the bounds of the spectrum have the same scaling with aspect ratio (sample size) as those for the whitened cross-correlation matrix~\cite{Bouchaud2007}. Thus, in this extremely well sampled limit, the cross-corelation and cross-covariance matrices can both be used to detect a signal.  However, the prefactor of this scaling is smaller for the cross-covariance matrix, indicating that whitening using the inverse of the empirically sampled self-covariance matrices introduces additional noise in the spectrum. Further, for sparse data, the cross-correlation cannot be evaluated---even if only one of the two variables is undersampled, where our results suggest that a signal may still be detectable in the cross-covariance. Together, these results suggest that in many cases the cross-covariance may be the most effective tool for detecting the shared signal in a pair of high-dimensional observations.

\begin{acknowledgments}
We thank  Philipp Fleig, Eslam Abdelaleem, and K.~Michael Martini for helpful discussions. This work was funded, in part, by a Simons Foundation Investigator grant, the NSF grant 2409416, and the NIH grant R01-NS084844.
\end{acknowledgments}

\appendix

\section{Calculating the spectrum of the empirical cross-covariance matrix}\label{append:sec1}

Here we calculate the spectrum of the $\nx \times \nx$-dimensional square of the normalized empirical cross-covariance matrix $\Csq$, given by Eq.~(\ref{eq:NECCM1}). 
Given $\nx ,\ \ny,\ T \gg 1$, this spectrum can be evaluated using random matrix theory. Parts of this calculation can be mapped onto previous  calculations~\cite{Philiplinear,Rocks2022,PhysRevE2010ref, dupic2014spectraldensityproductswishart, NIPS2017_Pennington} by reinterpreting the meaning of various variables. However, for pedagogical clarity, we choose to present a full, self-contained calculation here, which relies only on textbook RMT knowledge, instead of using special cases of calculations done with powerful, yet obscure mathematical machinery.

The nonzero eigenvalues of square of the NECCM $\Csq$ are the same as those of the matrix
\begin{align}\label{eq:Cprime_UV}
 \mathbf H&=\frac{1}{\sigma_{X}^2 \sigma_{\mathbf{Y}}^2 T^2} \left(\mathbf X\mathbf X^\top \right)\left(\mathbf Y\mathbf Y^\top\right)\nn \\
 &=\frac{\nx \ny}{ T^2} W_{X^\top} W_{Y^\top}\nn \\
 &=\frac{1}{p_X p_Y } W_{X^\top} W_{Y^\top}.
 \end{align}
Here $\mathbf{W_{X^\top}}$ and $\mathbf{W_{Y^\top}}$ are normalized Wishart matrices, given by
 \begin{align}\label{eq:white_Wishart}
     \mathbf W_{\mathbf X^\top}=\frac{1}{N_X \sigma_X^2}\mathbf X\mathbf X^\top\,,
\end{align}
and similar for $\mathbf{Y}$.  Crucially, $\mathbf{W_{X^\top}}$ and $\mathbf{W_{Y^\top}}$ are free matrices~\cite{freeness} (loosely, the appropriate generalization of statistical independence to noncommuting objects, such as  matrices). Freeness allows for certain matrix operations to commute with respect to each other. In classical probability, if $X$ and $Y$ are independent random variables, $\mathbb E[XY]=\mathbb E[X] \mathbb E[Y]$. Similarly, if $\mathbf{A}$ and $\mathbf{B}$ are free random matrices (in the large $N$ limit), then the limiting spectral distribution of the product $\mathbf{AB}$ or $\mathbf{A}^{1/2}\mathbf{B}\mathbf{A}^{1/2}$ can be obtained from the spectra of $\mathbf{A}$ and $\mathbf{B}$, in our case through Eqs.~(\ref{eq:S_multiplicative}, \ref{eq:S_scaling}).

The spectrum of $\matH$, $\rho_\matH$, can be evaluated from its Stieltjes transform, 
\begin{equation}
\mathfrak{h}(z)\equiv \mathfrak{g}_{\mathbf{H}}(z) \equiv \lim_{T\to \infty} \frac{1}{T}\mathbb E[\mathrm{Tr} (z\mathbf I- \mathbf H)],
\end{equation}
using the formula
\begin{align}\label{eq:sokhotski-plemelj1}
    \rho_\matH(\lambda)=\lim_{\eta \to 0^+}\frac{1}{\pi}\Im \mathfrak g_{\mathbf H}(\lambda - i\eta)\,.
\end{align}

To evaluate this Stieltjes transform, we must introduce the $\mathcal{T}$ and $\mathcal{S}$ transforms, which are useful for evaluating the Stieltjes transform of free products of random matrices (\cite{Potters_Bouchaud_2020}, Chapter.~15). The relevant properties of these transforms used in further calculations are summarized below.

The $\mathcal{T}$ transform of a matrix $A$ is defined as
 \begin{equation}\label{eq:T-transform_defn}
     \mathcal{T}_{\mathbf{A}}(z)=z\mathfrak g_{\mathbf{A}}(z)-1.
\end{equation}
The $\mathcal{T}$ transform, in turn, is used to define the $\mathcal S$ transform:
\begin{align}\label{eq:S_transformdefine}
    \mathcal S_\mathbf{A}(t)=\frac{t+1}{t{\mathcal T}_\mathbf{A}^{-1}(t)}.
\end{align}
For free matrices $\mathbf A$ and $\mathbf B$, the $\mathcal S$-transform of a product is multiplicative: \begin{align}\label{eq:S_multiplicative}
     \mathcal S_{\mathbf A\mathbf B}(t)=\mathcal S_{\mathbf A}(t)\mathcal S_{\mathbf B}(t)\,.
 \end{align}
Furthermore, for a scalar $a$,
\begin{align}\label{eq:S_scaling}
     \mathcal S_{a \mathbf A}(t)=a^{-1}\mathcal S_{\mathbf A}(t)\,.
 \end{align}

To derive the Stieltjes transform of $\matH$, we first evaluate its $\mathcal{S}$ transform. Using Eq.~(\ref{eq:S_multiplicative}) and Eq.~(\ref{eq:S_scaling}), we write
\begin{align}\label{eq:S_transformexpansioncase1}
   \mathcal{S}_{\matH}(t)&=\mathcal{S}\mleft(\frac{1}{p_X p_Y } \mathbf{W_{X^\top}} \mathbf{W_{Y^\top}}\mright)\nn\\
   &=p_X p_Y\mathcal{S}_{\mathbf{W_{X^\top}}}(t) \mathcal{S}_{\mathbf{W_{Y^\top}}}(t).
\end{align}
The $\mathcal{S}$-transform of a Wishart matrix is well known \cite{Potters_Bouchaud_2020}:
\begin{equation}\label{eq:StransformWishart}
 \mathcal{S}_{\mathbf{W_{X^\top}}}(t)=\frac{1}{1+p_Xt}.
\end{equation}

Now, plugging in the relevant terms for $\mathcal{S}_{\mathbf{W_{X}^\top}}(t)$ and $\mathcal{S}_{\mathbf{W_{Y}^\top}}(t)$ into Eq.~(\ref{eq:S_transformexpansioncase1})  and using Eq.~(\ref{eq:StransformWishart}), we obtain:
\begin{align}\label{eq:stransformexpansioncase1a}
     \mathcal{S}_{\matH}(t)=\frac{p_X}{1+p_Xt}\frac{p_Y}{1+p_Yt}.
\end{align}
To calculate the spectral density of the matrix of interest, we replace the $\mathcal S$-transform in Eq.~(\ref{eq:stransformexpansioncase1a}) with the corresponding $\mathcal{T}$-transform by using the relationship in Eq.~(\ref{eq:S_transformdefine}): 
\begin{align}\label{eq:inverseTcase1}
\mathcal{T}^{-1}_{\matH}(t)&=\frac{(t+1)(1+p_Xt)(1+p_Yt)}{tp_X p_Y}.
\end{align}
We now solve the equation for the functional inverse, ${\mathcal T}^{-1}(\mathcal T(z))=z$, using the definition of the $\mathcal T$-transform, Eq.~(\ref{eq:T-transform_defn}). This gives a cubic equation for the Stieltjes transform:
\begin{align}\label{eq:stilitsol}
   \mathfrak h^3&z^2p_Xp_Y+ \mathfrak h^2 z\left(p_Y(1-p_X)+p_X(1-p_Y)\right)\nn\\
   &+\mathfrak h \left((1-p_X)(1-p_Y)-z p_Xp_Y\right)\nn\\
   &+ p_Xp_Y=0\,.
 \end{align}

 Eq.~\ref{eq:stilitsol} can be obtained from the results in~\cite{dupic2014spectraldensityproductswishart,Rocks2022} by changing the definitions of parameters and rescaling variables appropriately. Ref~\cite{dupic2014spectraldensityproductswishart} further obtains the spectrum $\rho(\lambda)$ and studies its behavior in a few cases, but omits several important cases for data-science applications, such as the standard well-sampled case $T>\nx, \ny$ and the limiting behavior when the matrices have very different aspect ratios. 

The imaginary part of the roots of the cubic equation give us the density of eigenvalues. The bounds of the band $[\lambda_-, \lambda_+]$, for which the density is nonzero, are obtained from the zeros of the discriminant of the cubic equation. For an equation of the form
 \begin{equation}\label{eq:cubic}
    a \mathfrak h^3+  b \mathfrak h^2+ c \mathfrak h+  d=0,
 \end{equation}
the discriminant is
 \begin{equation}\label{eq:dis1}
     D=b^2c^2-4ac^3-4b^3d-27a^2d^2+18abcd,
 \end{equation}
where:
\begin{align}
    &a=z^2p_Xp_Y,\\
    &b=z\left(p_Y(1-p_X)+p_X(1-p_Y)\right),\\
    &c=\left((1-p_X)(1-p_Y)-z p_Xp_Y\right),\\
    &d=p_Xp_Y.
\end{align}

The density $\rho(\lambda)$ and the bounds $\lambda_{\pm}$ must then be transformed into the density of singular values $\rho(\gamma)$ and the bounds $\gamma_{\pm}$. For this, to get the spectrum of the nonzero part of the SVD of $\matC$, we use:
 \begin{equation}
     \rho_{A}(z)=2 z \rho_{A^2}(z^2),
 \end{equation}
 and the bounds obey $\gamma_{\pm}=\sqrt{\lambda_{\pm}}$. 
 
\section{The bounds of the spectrum in simplifying cases}

\subsection{Simplified solutions for $p_X=p_Y$ }

For $p_X=p_Y$, the cubic equation for the Stieltjes transform, Eq.~(\ref{eq:stilitsol}), reduces to:
\begin{align}\label{eq:stilitsolq_X=q_Y}
        &\mathfrak h^3z^2{p_X}^2+ \mathfrak h^2 z\left(p_X(1-p_X)+p_X(1-p_X)\right)\nn\\
   &+\mathfrak h \left((1-p_X)(1-p_X)-z {p_X}^2\right)+ p_X^2=0\,,
 \end{align}
and the discriminant (Eq.~\ref{eq:dis1}) simplifies to 
 \begin{align}
     D&=(4p_X^4- 12p_X^5+12p_X^6-4p_X^7)z^3\nn \\
     &+(-8p_X^6-20p_X^7+p_X^8)z^4+4p_X^8  z^5 \label{eq:xyequal_disc}.
 \end{align}

Solving Eq.~(\ref{eq:xyequal_disc}) for zeros we find that there are three zeros at $z=0$ and two zeroes at $z=z_{\pm}$, where  
\begin{equation}\label{eq:lamq_X=q_Y}
    z_{\pm}=\frac{8p_X^2 + 20 p_X^3 - p_X^4 \pm p_X^{5/2} (8 + p_X)^{3/2}}{8p_X^4}.
\end{equation}

When $T < N_X=N_Y$ ($p_X < 1$), the discriminant is negative in between $z_{\pm}$, and thus $\mathfrak{h}$ has a nonzero imaginary part, giving a nonzero density of eigenvalues. Thus, $\lambda_{\pm} = z_{\pm}$.
In the oversampled case, where $T > N_X=N_Y$ and thus $p_X > 1$, $z_-$ becomes negative. In this case, the discriminant is instead negative in the interval $(0, z_{+})$.
Thus, in general, we have
 \begin{align}
 \lambda_+ &= \frac{8p_X^2 + 20 p_X^3 - p_X^4 + p_X^{5/2} (8 + p_X)^{3/2}}{8p_X^4} \\ \lambda_{-} &= \begin{cases} \frac{8p_X^2 + 20 p_X^3 - p_X^4 + p_X^{5/2} (8 + p_X)^{3/2}}{8p_X^4}, & p_X < 1 \\ 0, & p_X \geq 1 . \end{cases} \label{eq:xyequal}
 \end{align}

For $p_X \gg 1$ (the comparison to whitened cross-correlation in the main text), we find
\begin{align}
    \lambda_+&\approx-\frac{1}{8}+\left( \frac{(8+p_X)^{3/2}}{8p_X^{3/2}}\right)\\
    &\approx\frac{3}{2p_X}.
\end{align}
Thus,
\begin{equation}
    \gamma_+=\sqrt{\frac{3}{2p_X}}.
\end{equation}

\subsection{Simplified solutions for $p_X < 1$, $p_Y \ll p_X$}\label{appendix:sec1b}

For $p_Y= \alpha p_X$ under the condition $\alpha \to 0$, the cubic equation for the Stieltjes transform Eq.~(\ref{eq:stilitsol}) reduces to:
\begin{align}\label{eq:stilitsolq_Y=epq_X}
  \alpha \mathfrak h^3&z^2p_X^2+ \mathfrak h^2 z p_X\left(\alpha(1-p_X)+(1-\alpha p_X)\right)\nn\\
   &+\mathfrak h \left((1-p_X)(1-\alpha p_X)-z \alpha p_X^2\right)+\alpha p_X^2=0\,.
 \end{align}
The discriminant of Eq.~(\ref{eq:stilitsolq_Y=epq_X}) is calculated using Eq.~(\ref{eq:dis1}). We then organize this discriminant as a polynomial in $z$, giving
 \begin{multline}\label{eq:dissol1}
 D=4 z^5 \alpha^4 p_X^8 +  z^4 (\alpha^2 p_X^6 + \alpha^3 (-10 p_X^6 - 10 p_X^7) \\
+ \alpha^4 (p_X^6 - 10 p_X^7 + p_X^8))+ 
 z^3 (\alpha (-2 p_X^4 - 2 p_X^5)\\
 + \alpha^2 (8 p_X^4 - 4 p_X^5 + 
    8 p_X^6) + \alpha^3 (-2 p_X^4 - 4 p_X^5 - 4 p_X^6 - 
    2 p_X^7)\\
+ \alpha^4 (-2 p_X^5 + 8 p_X^6 - 2 p_X^7) )
 +z^2 (p_X^2 - 2 p_X^3 + p_X^4 \\
 + \alpha (-2 p_X^2 + 2 p_X^3 + 2 p_X^4 - 
    2 p_X^5)  \\
+ \alpha^2 (p_X^2 + 2 p_X^3 - 6 p_X^4 + 2 p_X^5 + 
    p_X^6)\\
+ \alpha^3 (-2 p_X^3 + 2 p_X^4 + 2 p_X^5 - 
    2 p_X^6) + \alpha^4 (p_X^4 - 2 p_X^5 + 
    p_X^6)).
\end{multline}
Each term is of the form $f_n(\alpha)z^n$. As $\alpha \to 0$, we may expand each $f_n(\alpha)$ to the lowest nontrivial order in $\alpha$. Collecting the lowest-order terms for each power of $z$, the discriminant in Eq.~(\ref{eq:dissol1}) reduces to:
\begin{multline}\label{eq:dissol11}
 D\approx z^2 \left[p_X^2 (1-p_X)^2 -2(p_X^4+p_X^5) \alpha z + p_X^6 \alpha^2 z^2  \right. \\
 + \left.4  p_X^8 \alpha^4 z^3 \right].
\end{multline}

We seek positive roots $z_{\pm}(\alpha)$ of the right-hand group of terms (the equation has a single negative root, but since the eigenvalues of $\matH$ are positive by construction, this corresponds to a spurious root of the equation for $\mathfrak{h}$). This requires cancellation of at least two terms. That is, at least two terms of opposite signs must be of the same order in $\alpha$. We see that this can only happen if $z \sim \alpha^{-1}$ or $z\sim \alpha^{-3/2}$. In both of these possible cases, the final term is subleading and can be neglected. Thus, in this limit, we seek the roots of
\begin{align}\label{eq:dissol111}
 D&\approx(p_X^2 - 2 p_X^3 + p_X^4)z^2 -2(p_X^4+p_X^5)z^3 \alpha + z^4 p_X^6 \alpha^2.
\end{align}

We solve Eq.~(\ref{eq:dissol111}) for zeros. The \nth{4}-order equation has four zeroes. Two of the zeros are $z=0$, and the other two, $\lambda_{\pm}$, are
\begin{align}\label{eq:lambda1q_Y=epsilonq_X}
\lambda_{\pm}&=\frac{1+p_X\pm2\sqrt{p_X}}{\alpha p_X^2}\nn\\
&\approx\frac{1+p_X\pm2\sqrt{p_X}}{p_Xp_Y}.
\end{align}
Thus the  density of eigenvalues for for SVD of $\matC$ will be nonzero between $\gamma_{\pm}=\sqrt{\lambda_{\pm}}$, such that
\begin{align}\label{eq:gamma1q_Y=epsilonq_X}
\gamma_{\pm}&\approx\sqrt{\frac{1+p_X\pm2\sqrt{p_X}}{p_Xp_Y}}.
\end{align}

 \subsection{Simplified solutions for $p_X, p_Y \ll 1$}
 For $p_Y= \alpha p_X$ under the condition $p_X \to 0$ and $\alpha<1$, the cubic equation for the Stieltjes transform Eq.~(\ref{eq:stilitsol}) reduces to:
\begin{align}\label{eq:stilitsolq_Y=epq_Xq_X0}
  \alpha \mathfrak h^3&z^2p_X^2+ \mathfrak h^2 z p_X\left(\alpha+1\right)+\mathfrak h \left(1-z \alpha p_X^2\right)+\alpha p_X^2=0\,.
 \end{align}

The discriminant of  Eq.~(\ref{eq:stilitsolq_Y=epq_Xq_X0}) is calculated using Eq.~(\ref{eq:dis1}). Written as a polynomial in $z$, it is
\begin{multline}\label{eq:dissol2}
 D= 4 z^5 \alpha^4 p_X^8+ 
 z^4 (\alpha^2 p_X^6 - 
    10 \alpha^3 p_X^6 + \alpha^4 p_X^6 - 
    18 \alpha^3 p_X^7 \\
- 18 \alpha^4 p_X^7 - 
    27 \alpha^4 p_X^8) + 
 z^3 (-2 \alpha p_X^4 + 8 \alpha^2 p_X^4 - 
    2 \alpha^3 p_X^4 \\
- 4 \alpha p_X^5 + 
    6 \alpha^2 p_X^5 + 6 \alpha^3 p_X^5 - 
    4 \alpha^4 p_X^5) \\
+ z^2 (p_X^2 - 2 \alpha p_X^2 + \alpha^2 p_X^2).
\end{multline}
As $p_X \to 0$, the contribution of higher-order terms for each power of $z$ to the final solution will be negligible. Collecting the lowest order terms in $p_X$ for each power of $z$, the discriminant in Eq.~(\ref{eq:dissol2}) reduces to
\begin{multline}\label{eq:dissol21}
 D= 4 z^5 \alpha^4 p_X^8+ 
 z^4 (\alpha^2 p_X^6 - 
    10 \alpha^3 p_X^6 + \alpha^4 p_X^6) \\
+ z^3 (-2 \alpha p_X^4 + 8 \alpha^2 p_X^4 - 
    2 \alpha^3 p_X^4)+\\ z^2 (p_X^2 - 2 \alpha p_X^2 + \alpha^2 p_X^2).
\end{multline}
We solve Eq.~(\ref{eq:dissol21}) for zeros. The \nth{5}-order equation has 5 zeroes (counting their multiplicities). Two of the zeroes are at $z=0$, one is at $z=\frac{-(1-\alpha)^2}{4 \alpha^2 p_X^2} < 0$. Thus, the other two are $\lambda_{\pm}$. Taking the condition $D<0$, we find that nonzero density requires $\lambda \in [\lambda_-, \lambda_+]$. In particular, we find the solution
\begin{equation}
\lambda_{\pm}=\frac{1\pm2\sqrt{p_X(1+\alpha)}}{\alpha p_X^2}.
\end{equation}
The nonzero density of eigenvalues for  SVD of $\matC$ will be between $\gamma_{\pm}=\sqrt{\lambda_{\pm}}$, where
\begin{align}  \gamma_{\pm}=\sqrt{\lambda_{\pm}}&=\sqrt{\frac{1\pm2\sqrt{p_X(1+\alpha)}}{\alpha p_X^2}}\nn\\
    &=\sqrt{\frac{1\pm2\sqrt{p_X+p_Y}}{p_Y p_X}}\nn\\
    &\approx\frac{1\pm\sqrt{p_Y+p_X}}{\sqrt{p_Yp_X}}.\label{eq:pxpy<<1}
\end{align}
In the final step we have again used the fact that we are studying the special case $p_X, p_Y \ll 1$, and $\sqrt{1+x} = 1 + x/2 + O(x^2)$.

Setting $p_X=p_Y$ in Eq.~(\ref{eq:pxpy<<1}) gives Eq.~(\ref{eq:edge1}) in the main text.

\bibliography{main}

\end{document}